\documentclass{amsart}

\usepackage[latin1]{inputenc}
\usepackage{amsmath}
\usepackage{amssymb}
\usepackage[pdftex]{graphicx}
\usepackage{epsfig}
\usepackage{url}
\usepackage[breaklinks=true,hyperref]{hyperref}

\newtheorem{theorem}{Theorem}[section]
\newtheorem{proposition}[theorem]{Proposition}
\newtheorem{corollary}[theorem]{Corollary}
\newtheorem{definition}[theorem]{Definition}
\newtheorem{example}[theorem]{Example}
\newtheorem{conjecture}[theorem]{Conjecture}
\newtheorem{remark}[theorem]{Remark}
\newtheorem{openproblem}[theorem]{Open Problem}

\begin{document}
	\title{$q,t$-Fu\ss-Catalan numbers for complex reflection groups}
	\author{Christian Stump}
	\address{Fakult\"at f\"ur Mathematik, Universit\"at Wien, Nordbergstra{\ss}e 15, A-1090 Vienna, Austria}
	\email{christian.stump@univie.ac.at}
	\urladdr{http://homepage.univie.ac.at/christian.stump/}
	\date{\today}
	\thanks{Research supported by the Austrian Science Foundation FWF, grant P17563-N13 "Macdonald polynomials and q-hypergeometric Series"}
	
	\begin{abstract}
		In type $A$, the $q,t$-Fu\ss -Catalan numbers $\operatorname{Cat}_n^{(m)}(q,t)$ can be defined as a bigraded Hilbert series of a module associated to the symmetric group $\mathcal{S}_n$. We generalize this construction to (finite) complex reflection groups and exhibit some nice conjectured algebraic and combinatorial properties of these polynomials in $q$ and $t$. Finally, we present an idea how these polynomials could be related to some graded Hilbert series of modules arising in the context of rational Cherednik algebras. This is work in progress.\\ \\
		\noindent {\sc R\'esum\'e.}
		Dans le cas du type $A$, les $q,t$-nombres de Fu\ss -Catalan $\operatorname{Cat}_n^{(m)}(q,t)$ peuvent \^etre d\'efinis comme la s\'erie de Hilbert bigradu\'ee d'un certain module associ\'e au groupe sym\'etrique $\mathcal{S}_n$. Nous g\'en\'eralisons cette construction aux groupes de r\'eflexion complexes (finis) et nous formulons de jolies propri\'et\'es (conjecturales) alg\'ebriques et combinatoires de ces polyn\^omes en $q$ et $t$. Enfin, nous d\'ecrivons une id\'ee sur la mani\`ere dont ces polyn\^omes pourraient \^etre li\'es \`a certaines s\'eries de Hilbert de modules apparaissant dans le contexte des alg\`ebres de Cherednik rationnelles. Ceci est un travail en cours.
	\end{abstract}

	\maketitle

\section{Introduction}
		Within the last 15 years the $q,t$-Fu\ss -Catalan numbers of type $A$, $\operatorname{Cat}_n^{(m)}(q,t)$, arose in more and more contexts in mathematics, namely in \emph{symmetric functions theory}, \emph{algebraic} and \emph{enumerative combinatorics}, \emph{representation theory} and \emph{algebraic geometry}. They first appeared in a paper by A.~Garsia and M.~Haiman in the context of \emph{modified Macdonald polynomials}, \cite{garsiahaiman}. Later, in his work on the $n!$-- and on the $(n+1)^{n-1}$~--conjecture, M.~Haiman showed that $\operatorname{Cat}_n^{(m)}(q,t)$ is equal to the Hilbert series of the alternating component of the \emph{diagonal coinvariant ring}, \cite{haiman3}. J.~Haglund, \cite{haglund2}, and N.~Loehr, \cite{loehr}, found a very interesting (partially conjectured) combinatorial interpretation of $\operatorname{Cat}_n^{(m)}(q,t)$. The $q,t$-Catalan numbers have many interesting algebraic and combinatorial properties. To mention some: they are symmetric functions in $q$ and $t$ with positive integer coefficients and give a $q,t$-extension of the famous \emph{Fu\ss -Catalan numbers}
		$$\operatorname{Cat}_n^{(m)} := \frac{1}{mn+1}{\binom{(m+1)n}{n}}.$$
		Specializing $t=1$ in $\operatorname{Cat}_n^{(m)}(q,t)$ gives some $q$-Fu\ss -Catalan numbers introduced by J.~F\"urlinger and J.~Hofbauer in \cite{fuerlingerhofbauer} and specializing $t=q^{-1}$ gives another $q$-extension of the Fu\ss -Catalan numbers introduced in \cite{macmahon} by P.A.~MacMahon.\\ \\
		We generalize the definition of Fu\ss-Catalan numbers by means of a Hilbert series from the symmetric group to arbitrary finite complex reflection groups. Furthermore, we present conjectured generalizations of the properties mentioned above concerning specializations of $q$ and $t$, and combinatorial interpretations, see Conjectures \ref{conjecture1}, \ref{conjecture2} and \ref{conjecture4}. Finally, we present an idea which would relate them to some graded Hilbert series constructed by I.~Gordon \cite{gordon} and by Y.~Berest, P.~Etingof and V.~Ginzburg \cite{BEG} in the context of \emph{rational Cherednik algebras}, see Conjecture \ref{conjecture3} and the following corollary.\\ \\
		This extended abstract is organized as follows: in Section \ref{sec1}, we define alternating polynomials of type $A$ and generalize this definition to complex reflection groups. In Section \ref{sec2}, we first define $q,t$-Fu\ss -Catalan numbers of type $A$ as a bigraded Hilbert series and present some properties. Then we generalize this definition to $q,t$-Fu\ss -Catalan numbers for complex reflection groups, $\operatorname{Cat}_n^{(m)}(W,q,t)$, and present two conjectures concerning the specializations $q=t=1$ and $t=q^{-1}$. In Section \ref{sec3}, we combinatorially define $q$-Fu\ss -Catalan numbers in terms of the \emph{extended Shi arrangement} and conjecture that they appear in the specialization $t=1$ of $\operatorname{Cat}_n^{(m)}(W,q,t)$. In Section \ref{sec4}, we present another conjecture which would connect $\operatorname{Cat}_n^{(m)}(W,q,t)$ to a graded Hilbert series in the context of \emph{rational Cherednik algebras}.
	\section{Alternating polynomials}\label{sec1}
	\subsection{Alternating polynomials associated to the symmetric group} \label{sec11}
		The symmetric group $\mathcal{S}_n$, which is the reflection group of type $A_{n-1}$, acts on the polynomial ring
		$$\mathbb{C}[\mathbf{x},\mathbf{y}] := \mathbb{C}[x_1,y_1, \ldots ,x_n,y_n]$$
		by permuting the coordinates in $\mathbf{x}$ and $\mathbf{y}$ simultaneously amongst themselves. This is the \emph{diagonal action}
		$$\sigma(x_i):=x_{\sigma(i)},\sigma(y_i):=y_{\sigma(i)} \hspace{10pt} \text{for } \sigma \in \mathcal{S}_n.$$
		A polynomial $p \in \mathbb{C}[\mathbf{x},\mathbf{y}]$ is called \emph{alternating} if
		$$\sigma(p) = \operatorname{sgn}(\sigma) p \hspace{10pt} \text{for all } \sigma \in \mathcal{S}_n,$$
		where $\operatorname{sgn}(\sigma)$ is the usual sign of the permutation $\sigma$. We denote the space of all alternating polynomials by $\mathbb{C}[\mathbf{x},\mathbf{y}]^\epsilon$.
	As a vector space, $\mathbb{C}[\mathbf{x},\mathbf{y}]^\epsilon$ has a well-known basis: For $G=\{(\alpha_1,\beta_1),\ldots,(\alpha_n,\beta_n)\} \subseteq \mathbb{N} \times \mathbb{N}$ define a bivariate analogue of the Vandermonde determinant $\Delta_G$ by
	$$\Delta_G=\det \left( \begin{array}{ccc} x_1^{\alpha_1}y_1^{\beta_1} & \ldots & x_1^{\alpha_n}y_1^{\beta_n} \\
	\vdots & & \vdots \\
	x_n^{\alpha_1}y_n^{\beta_1} & \ldots & x_n^{\alpha_n}y_n^{\beta_n} \end{array} \right).$$
	The set
	$$\mathcal{B} := \big\{\Delta_G : G \subseteq \mathbb{N} \times \mathbb{N}, |G|=n \big\}$$
	forms a vector space basis of $\mathbb{C}[\mathbf{x},\mathbf{y}]^\epsilon$ and in particular the ideal generated by $\mathcal{B}$ is the same as the ideal generated by all alternating polynomials.
		
	\subsection{Alternating polynomials associated to any complex reflection group}\label{refrep}
	The concept for polynomials to be alternating can be generalized to any (finite) complex reflection group in the following way: let $V$ be an $n$-dimensional complex vector space and let $W \subseteq \operatorname{GL}(V)$ be a (finite) complex reflection group acting on $V$. For definitions and further information on complex reflection groups see e.g. \cite{bessis}.\\ \\
	The \emph{contragredient action} of $W$ on $V^\ast = \operatorname{Hom}(V,\mathbb{C})$ is given by
	$$\omega(\rho):= \rho \circ \omega^{-1}.$$
	This induces an action of $W$ on the symmetric algebra $S(V^\ast)$ which is equal to $\mathbb{C}[\mathbf{x}]$. \lq\lq Doubling up\rq\rq\ this action diagonally defines a \emph{diagonal action} of $W$ on $\mathbb{C}[\mathbf{x},\mathbf{y}]$.
	\begin{definition}
		Let $W$ be a complex reflection group acting on a complex vector space of dimension $n$. We call a polynomial $p \in \mathbb{C}[\mathbf{x},\mathbf{y}]$ \emph{alternating} if
		$$ \det(\omega) \omega(p) = p \hspace{10pt} \text{for all } \omega \in W.$$
	\end{definition}
	For $W$ being the complex reflection group of type $A_{n-1}$ - which is the symmetric group $\mathcal{S}_n$ - this definition reduces to the definition of alternating polynomials given above. We denote the space of all alternating polynomials by $\mathbb{C}[\mathbf{x},\mathbf{y}]^{\epsilon,W}$.
	\begin{remark}
		If $W$ is a real reflection group or, equivalently, if $W$ is a (finite) Coxeter group then
		$$\det(\omega) = (-1)^{l(\omega)},$$
		where $l$ is the length function in the Coxeter group $W$.
	\end{remark}
	The reason why we call a polynomial $p$ alternating if $\det(\omega)  \omega(p) = p$ and not if $\omega(p) = \det(\omega) p$ is the following: Define the \emph{sign idempotent} $\mathbf{e}_\epsilon$ by
	$$\mathbf{e}_\epsilon := \frac{1}{|W|}\sum_{\omega \in W}{\det(\omega) \omega}$$
	and the \emph{sign representation} $\epsilon$ by $ \omega(z) := \det(\omega) z$ for all $\omega \in W$ and $z \in \mathbb{C}$. Then
	$$\mathbb{C}[\mathbf{x},\mathbf{y}]^{\epsilon,W} = (\mathbb{C}[\mathbf{x},\mathbf{y}] \otimes \epsilon)^W = \mathbf{e}_\epsilon\mathbb{C}[\mathbf{x},\mathbf{y}].$$
	As for the symmetric group, $\mathbb{C}[\mathbf{x},\mathbf{y}]^{\epsilon,W}$ has a vector space basis given by
	$$\mathcal{B}_W := \big\{\mathbf{e}_\epsilon(m) : m \text{ monomial in } \mathbf{x},\mathbf{y} \text{ with } \mathbf{e}_\epsilon(m) \neq 0 \big\}.$$
	\begin{remark}
		\begin{itemize}
			\item For $W$ of type $A$, $\mathcal{B}_W$ reduces to $\mathcal{B}$ defined above,
			\item for $W$ of type $B$, $\mathcal{B}_W$ reduces to $\big\{\Delta_G : G \subseteq \mathbb{N} \times \mathbb{N}, |G|=n, \alpha_i+\beta_i \equiv 1 \operatorname{mod} 2 \big\}$,
			\item for $W$ of type $D$, $\mathcal{B}_W$ reduces to $\big\{\Delta_G : G \subseteq \mathbb{N} \times \mathbb{N}, |G|=n, \alpha_i+\beta_i \equiv \alpha_j+\beta_j \operatorname{mod} 2 \big\}$.
		\end{itemize}
	\end{remark}
	
	\section{$q,t$-Fu\ss -Catalan numbers}\label{sec2}
	Before we define $q,t$-Fu\ss -Catalan numbers in general, we review the definition and the properties about the well-studied case $W=\mathcal{S}_n$ which they seem to generalize. To refer to the parameter $m \in \mathbb{N}$, we use the term \emph{Fu\ss -Catalan} which has commonly been used in the literature for \emph{higher} Catalan numbers of general type, see e.g. \cite{armstrong} and \cite{fomin}. In the literature concerning only the case $W=\mathcal{S}_n$, the name \emph{generalized} $q,t$\emph{-Catalan numbers} was more usual.
	\subsection{$q,t$-Fu\ss -Catalan numbers associated to the symmetric group}
	Let $\mathcal{S}_n$ act on $\mathbb{C}[\mathbf{x},\mathbf{y}]$ as described in Section \ref{sec11} and let $I \trianglelefteq \hspace{2pt}  \mathbb{C}[\mathbf{x},\mathbf{y}]$ be the ideal generated by all alternating polynomials. Define the $\mathcal{S}_n$-module $M^{(m)}$ to be the \emph{minimal generating space} of $I^m$,
		$$M^{(m)}:=I^m/\langle \mathbf{x},\mathbf{y} \rangle I^m.$$
		It carries a natural bigrading by degree in $\mathbf{x}$ and degree in $\mathbf{y}$, $M^{(m)} = \oplus_{i,j \geq 0}{M^{(m)}_{ij}}$.
		\begin{remark}
			The name \emph{minimal generating space} comes from the fact that - as a vector space - $M^{(m)}$ is isomorphic to the complex vector space with basis in one-to-one correspondence to any homogeneous minimal generating set of $I^m$.
		\end{remark}
		The following definition is due to M.~Haiman:
		\begin{definition}
			The $q,t$\emph{-Fu\ss -Catalan numbers of type} $A_{n-1}$ are defined as the bigraded Hilbert series of the $\mathcal{S}_n$-module $M^{(m)}$,
			$$\operatorname{Cat}^{(m)}_n(q,t) := \mathcal{H}{(M^{(m)};q,t)} = \sum_{i,j \geq 0}{\dim{(M^{(m)}_{ij})}q^it^j}.$$
		\end{definition}
	\begin{remark}
		$\operatorname{Cat}^{(m)}_n(q,t)$ was originally defined in \cite{garsiahaiman} as the complicated rational function
		$$\sum_{\mu \vdash n}{\frac{q^{(m+1)n(\mu')}t^{(m+1)n(\mu)}(1-q)(1-t)\Pi_{\mu}{(q,t)}B_\mu(q,t)}{\prod_{c \in D(\mu)}{(q^{a(c)}-t^{l(c)+1})(t^{l(c)}-q^{a(c)+1})}}}$$
		in the context of \emph{modified Macdonald polynomials}. M.~Haiman later showed that this rational function is in fact equal to the bigraded Hilbert series in the definition above, \cite{haiman5}.
	\end{remark}
	In \cite{garsiahaglund2}, A.~Garsia and J.~Haglund proved a simple combinatorial interpretation of $\operatorname{Cat}_n(q,t) := \operatorname{Cat}^{(1)}_n(q,t)$ which was conjectured by J.~Haglund in \cite{haglund2}, where he introduced the \emph{bounce statistic} on \emph{Catalan paths} or, equivalently, on the set $\mathcal{D}_n$ of all  partitions that fit inside the partition $(n-1,\ldots,2,1)$:
	$$\operatorname{Cat}_n(q,t) = \sum_{\lambda \in \mathcal{D}_n}{q^{\operatorname{area}(\lambda)}t^{\operatorname{bounce}(\lambda)}}.$$
	Together with N.~Loehr, they extended the definitions of \emph{area} and \emph{bounce} to the set $\mathcal{D}^{(m)}_n$ of $m$-Catalan paths, these are partitions that fit inside the partition $((n-1)m,\ldots,2m,m)$, and conjectured a combinatorial interpretation of $\operatorname{Cat}^{(m)}_n(q,t)$ in terms of these statistics, \cite{loehr}:
	$$\operatorname{Cat}_n^{(m)}(q,t) = \sum_{\lambda \in \mathcal{D}_n^{(m)}}{q^{\operatorname{area}(\lambda)}t^{\operatorname{bounce}(\lambda)}}.$$
	\begin{remark}
			The specialization $t=1$ was proved by A.~Garsia and M.~Haiman in \cite{garsiahaiman}.
	\end{remark}
	\begin{corollary}
		The specialization $q=t=1$ reduces $\operatorname{Cat}_n^{(m)}(q,t)$ to the $\operatorname{Cat}_n^{(m)}$.
	\end{corollary}
	Also in \cite{garsiahaiman}, it was shown that the specialization $t=q^{-1}$ yields the following $q$-extension of $\operatorname{Cat}_n^{(m)}$:
	\begin{theorem}[Garsia, Haiman]\label{theo1}
		$$q^{m\binom{n}{2}}\operatorname{Cat}_n^{(m)}(q,q^{-1}) = \frac{1}{[mn+1]_q}\begin{bmatrix} (m+1)n \\ n \end{bmatrix}_q,$$
		where $[k]_q := 1+q+\ldots+q^{k-1}, [k]_q!:= [1]_q[2]_q \cdots [k]_q$ and $\left[\begin{smallmatrix} k \\ l \end{smallmatrix}\right]_q:=[k]_q!/[l]_q![k-l]_q!$ .
	\end{theorem}
	\begin{remark}
		For $m=1$, this reduces to the \emph{q-Catalan numbers} defined by P.A.~MacMahon in \cite{macmahon}.
	\end{remark}
	
	\subsection{$q,t$-Fu\ss -Catalan numbers for complex reflection groups}\label{defqtfusscat}
	Recall that a complex reflection group $W$ acts on the polynomial ring $\mathbb{C}[\mathbf{x},\mathbf{y}]$ diagonally as described in Section \ref{refrep}. Let $I \trianglelefteq \hspace{2pt}  \mathbb{C}[\mathbf{x},\mathbf{y}]$ be the ideal generated by all alternating polynomials and define the $W$-module $M^{(m)} := I^m/ \langle \mathbf{x},\mathbf{y} \rangle I^m$.
	\begin{definition}
			The $q,t$\emph{-Fu\ss -Catalan numbers} associated to $W$ are defined as
			$$\operatorname{Cat}^{(m)}(W,q,t) := \mathcal{H}{(M^{(m)};q,t)} = \sum_{i,j \geq 0}{\dim{(M_{ij})}q^it^j}.$$
	\end{definition}
	\begin{remark}
		Using the computer algebra systems \emph{Singular} and \emph{Macaulay 2}, we computed the dimension of $M^{(m)}$ as well as  $\operatorname{Cat}^{(m)}(W,q,t)$ for the classical types at least up to rank $4$ and small $m$ and following exceptional types
		$$I_2(k) \text{ for } k \in \{5,6,10,12\},\hspace{5pt} H_3,\hspace{5pt} G(k,1,1) \hspace{5pt} \text{for } k \leq 10,\hspace{5pt} G(4,2,2).$$
		All following conjectures are based on these computations.
	\end{remark}
	In \cite{reiner2}, V. Reiner defined \emph{Fu\ss -Catalan numbers} for classical reflection groups, and in \cite{bessis2} D. Bessis generalized this definition to well-generated complex reflection group $W$:
	$$\operatorname{Cat}^{(m)}(W) := \prod_{i=1}^l{\frac{d_i+mh}{d_i}},$$
	where $l$ is the \emph{rank} of $W$, $d_1 \leq \ldots \leq d_l$ are the \emph{degrees of the fundamental invariants} and $h = d_l$ is the \emph{Coxeter number}. For definitions and further information see e.g. \cite[Section 2.7]{armstrong} and \cite{bessis2}.\\ \\
		For more or less general classes of reflection groups, this number counts a bunch of interesting combinatorial objects, see e.g. \cite{armstrong} and \cite{atha}. For $W=\mathcal{S}_n$, it reduces to $\operatorname{Cat}^{(m)}_n$ and for any real reflection group, it reduces for $m=1$ to the well-known Catalan numbers $\operatorname{Cat}(W)$ for real reflection groups, which are shown in Fig.~\ref{Cat}.
	\begin{figure}
		\centering
		\begin{tabular}[h]{|c|c|c|c|c|c|c|c|c|c|}
			\hline
			$A_{n-1}$ & $B_n$ & $D_n$ & $I_2(k)$ & $H_3$ & $H_4$ & $F_4$ & $E_6$ & $E_7$ & $E_8$ \\
			\hline
			$\frac{1}{n+1}{\binom{2n}{n}}$ & ${\binom{2n}{n}}$ & ${\binom{2n}{n}} - \binom{2(n-1)}{n-1}$ & $k+2$ & $32$ & $280$ & $105$ & $833$ & $4160$ & $25080$ \\
			\hline
		\end{tabular}
		\caption{$\operatorname{Cat}(W)$ for the irreducible real reflection groups.}
		\label{Cat}
	\end{figure}
	\newpage
	Our first conjecture concerns $\operatorname{Cat}^{(m)}(W,1,1)$, the dimension of $M$: 
	\begin{conjecture}\label{conjecture1}
		Let $W$ be a well-generated complex reflection group. Then
		$$\operatorname{Cat}^{(m)}(W,1,1) = \operatorname{Cat}^{(m)}(W).$$
	\end{conjecture}
	In \cite{bessis}, D. Bessis and V. Reiner defined a $q$-extension of $\operatorname{Cat}^{(m)}(W)$ by
	$$\prod_{i=1}^l{\frac{[d_i+mh]_q}{[d_i]_q}}.$$
	The following conjecture, which is obviously stronger than Conjecture \ref{conjecture1}, would generalize Theorem \ref{theo1} and would thereby give a new answer to a question of C. Kriloff and V. Reiner in \cite[Problem 2.2]{reiner}:
	\begin{conjecture}\label{conjecture2}
		Let $W$ be a well-generated complex reflection group. Then
	$$q^{mN} \operatorname{Cat}^{(m)}(W,q,q^{-1}) = \prod_{i=1}^l{\frac{[d_i+mh]_q}{[d_i]_q}},$$
	where $N := \sum{(d_i-1)}$.
	\end{conjecture}
	\begin{remark}
		For a real reflection group $W$, the $(d_i-1)$'s appearing in the conjecture are the \emph{exponents} associated to $W$ and $N$ is equal to the number of positive roots.
	\end{remark}
	\begin{openproblem}\label{problem1}
		Are there statistics \emph{qstat} and \emph{tstat} on objects counted by $\operatorname{Cat}^{(m)}(W)$ which generalize \emph{area} and \emph{bounce} on Catalan paths $\mathcal{D}_n^{(m)}$ such that
	$$\operatorname{Cat}^{(m)}(W,q,t) = \sum_{\lambda}{q^{\operatorname{qstat}(\lambda)} t^{\operatorname{tstat}(\lambda)}}?$$
	\end{openproblem}
	In the next section, we will present some conjectures concerning this open problem.
\section{A generalization of the \emph{area} statistic to crystallographic reflection groups}\label{sec3}
\subsection{The \emph{area} statistic for Catalan paths of type $A$} \label{CatpathsA}
	Specializing $t=1$ reduces $\operatorname{Cat}^{(m)}_n(q,t)$ to the well-known \emph{Carlitz} $q$\emph{-Fu\ss-Catalan numbers} defined by
	$$\operatorname{Cat}_n^{(m)}(q) := \operatorname{Cat}_n^{(m)}(q,1) = \sum_{\lambda \in \mathcal{D}_n^{(m)}}{q^{\operatorname{area}(\lambda)}}.$$
	They satisfy the following recurrence which can be deduced from a generating function identity proved by C. Krattenthaler in \cite[Theorem 9]{kratt1}:
	\begin{theorem}
		$$\operatorname{Cat}_{n+1}^{(m)}(q) = \sum_{k_1 + \ldots + k_{m+1} = n}{q^{n(\mathbf{k})} \operatorname{Cat}_{k_1}^{(m)}(q) \ldots \operatorname{Cat}_{k_{m+1}(q)}^{(m)}}, \qquad \operatorname{Cat}_{0}^{(m)}(q) = 1,$$
		where $n(\mathbf{k}) = n(k_1,\ldots,k_{m+1}) := \sum{(m+1-i) k_i}$.
	\end{theorem}

\subsection{An \emph{area} statistic for Catalan paths of type $B$}  \label{CatpathsB}
\begin{figure}
		\centering
		\includegraphics[width=350pt]{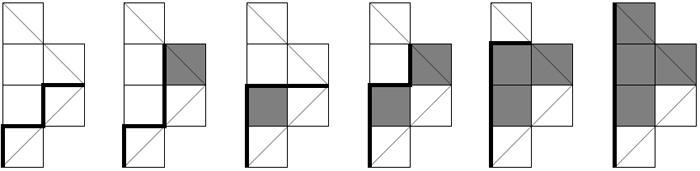}
		\caption{All type $B$ Catalan paths of length $2$.}
		\label{Bpaths}
	\end{figure}

	For $m=1$, we define an \emph{area}-statistic on type $B$ Catalan paths and establish an analogous recurrence.
	\begin{definition}
		A \emph{type} $B$ \emph{Catalan path of length} $n$ is a lattice paths of $2n$ steps, either north or east, that starts at $(0,0)$ and stays above the diagonal $x=y$. For such a path $\lambda$, we define $\operatorname{area}(\lambda)$ to be the number of boxes in the region confined by the path, the diagonal $x=y$ and the anti-diagonal $x=2n-y$, not counting the halfboxes at the diagonal $x=y$ but counting the halfboxes at the anti-diagonal $x=2n-y$.
	\end{definition}
	\begin{example}\label{ex1}
		In Fig.~\ref{Bpaths}, all Catalan paths of type $B_2$ are shown, the boxes which contribute to the area are shaded.
	\end{example}
	In analogy to type $A$, we define $q$-\emph{Catalan numbers} for type $B$ in the following way:
	\begin{definition}
		$$\operatorname{Cat}_{B_n}(q) := \sum{q^{\operatorname{area}(\lambda)}},$$
		where the sum ranges over all type $B$ Catalan paths $\lambda$ of length $n$.
	\end{definition}
	\begin{example}
		As shown in Example \ref{ex1}, we have $\operatorname{Cat}_{B_2}(q) = 1 + 2q+q^2+q^3+q^4.$
	\end{example}
	The definition is based on the following conjecture:
	\begin{conjecture}
		$$\operatorname{Cat}^{(1)}(W_{B_n},q,1) = \operatorname{Cat}_{B_n}(q).$$
	\end{conjecture}
	$\operatorname{Cat}_{B_n}(q)$ satisfy the following recurrence involving Catalan numbers of type $A$:
	\begin{theorem}\label{theo2}
		$$\operatorname{Cat}_{B_n}(q) = \operatorname{Cat}_n(q) + \sum_{k=0}^{n-1}{q^{2k+1} \operatorname{Cat}_{B_k}(q) \operatorname{Cat}_{n-k}(q)}, \qquad \operatorname{Cat}_{B_0}(q) = 1.$$
	\end{theorem}
	\begin{proof}
		Let $\lambda$ be a type $B$ Catalan path of length $n$. Then either $\lambda$ has as many east as north steps, which means $\lambda$ is equal to a type $A$ Catalan path of length $n$, or there exists a last point $(k,k+1)$ where the path touches the diagonal $x+1=y$ and stays strictly above afterwards. Now, we have an initial type $A$ like Catalan path of length $k+1$ (where the last step is a north step instead of an east step). After this north step, a type $B$ Catalan path of length $n-k-1$ starts. This gives the following recurrence which is equivalent to the statement:
		$$\operatorname{Cat}_{B_n}(q) = \operatorname{Cat}_n(q) + \sum_{k=0}^{n-1}{q \operatorname{Cat}_{k+1}(q) q^{2(n-k-1)}\operatorname{Cat}_{B_{n-k-1}}(q)}.$$
	\end{proof}
	\begin{corollary}
		The $\operatorname{Cat}_{B_n}(q)$ satisfy the following generating function identity:
		$$\sum_{n \geq 0}{\frac{x^n q^{-n(n-1)}(1-qx)}{(-x;q^{-1})_{2n+1}}\operatorname{Cat}_{B_n}(q)} = 1,$$
		where $(a;q)_k := (1-a)(1-qa)\ldots (1-q^{k-1}a)$.
	\end{corollary}
		We will see in the next section that both, $m$-Catalan paths of type $A$ and Catalan paths of type $B$ are special cases of a more general construction and that it is not possible to construct Catalan paths of type $B$ for higher $m$'s as lattice paths (at least not in the manner of defining an area generating function equal to the specialization $t=1$ of $q,t$-Fu\ss -Catalan numbers). This is likely to be the reason why we were - so far - not able to find a recurrence in type $B$ for higher $m$'s.
	\subsection{The extended Shi arrangement and the coheight statistic}
	Fix $\Phi$ to be a crystallographic root system and let $W=W_\Phi$ be the associated reflection group. The \emph{root poset} of $\Phi$ is given by the partial order on the set of positive roots $\Phi^+$ defined by covering relation
		$$\alpha \prec \beta :\Leftrightarrow \beta - \alpha \text{ is a simple root}.$$
	An \emph{order ideal} $I \trianglelefteq \hspace{2pt} \Phi^+$ is a subset $I \subseteq \Phi^+$ such that $\alpha \leq \beta \in I \text{ implies } \alpha \in I.$
	\begin{theorem}[V. Reiner \cite{reiner2}] Let $\Phi$ be a crystallographic root system. Then
	$$\#\{I \trianglelefteq \hspace{2pt} \Phi^+\} = \operatorname{Cat}(W).$$
	\end{theorem}
	In \cite{atha}, C.A. Athanasiadis generalized this theorem to $\operatorname{Cat}^{(m)}(W)$ as follows:
	Let $\mathcal{I}$ be an increasing chain of order ideals $I_1 \subseteq \ldots \subseteq I_m \subseteq \Phi^+$. We call $\mathcal{I}$ a \emph{filtered chain} of length $m$ if for $i,j \geq 1$,
	\begin{eqnarray*}
		(I_i + I_j) \cap \Phi^+ &\subseteq& I_{i+j} \hspace{10pt}\text{with } i+j \leq m,\\
		(J_i + J_j) \cap \Phi^+ &\subseteq& J_{i+j},
	\end{eqnarray*}
	where $J_i = \Phi^+ \setminus I_i$ and $J_i = J_m$ for $i > m$.\\ \\
	Let $V$ be the vector space spanned by $\Phi$, with inner product $($ , $)$. The \emph{extended Shi arrangement} $\mathsf{Shi}^{(m)}(\Phi)$ is given by the collection of hyperplanes in $V$ defined by the affine equations $(\alpha,x) = k$ for $\alpha \in \Phi$ and $-m < k \leq m$. Thus $\mathsf{Shi}^{(m)}(\Phi)$ is a deformation of the Coxeter arrangement $\mathcal{A}_\Phi$. A \emph{positive region of} $\mathsf{Shi}^{(m)}(\Phi)$ is a connected component of $V \setminus \mathsf{Shi}^{(m)}(\Phi)$ which lies in the fundamental chamber of $\mathcal{A}_\Phi$.\\ \\
	C.A. Athanasiadis defined the following map $\Psi$ between positive regions of $\mathsf{Shi}^{(m)}(\Phi)$ and filtered chains in $\Phi^+$ and showed that $\Psi$ is a bijection: let $R$ be a positive region and let $x \in R$. Then $\Psi(R)$ is defined to be the filtered chain $I_1 \subseteq \ldots \subseteq I_m \subseteq \Phi^+$ such that
	\begin{eqnarray*}
		(\alpha,x) < i &,& \text{if } \alpha \in I_i \\
		(\alpha,x) > i &,& \text{if } \alpha \in J_i = \Phi^+ \setminus I_i.
	\end{eqnarray*}
	\begin{theorem}[C.A. Athanasiadis \cite{atha}] Both the number of filtered chains in $\Phi^+$ and the number of positive regions of $\mathsf{Shi}^{(m)}(\Phi)$ is equal to $\operatorname{Cat}^{(m)}(W)$.
	\end{theorem}
	\begin{definition}
		Let $R^0$ be the \emph{fundamental region} given by $0 < (x,\alpha) < 1$ for all $\alpha \in \Phi^+$. For any region $R$, define the \emph{height} of $R$, denoted $\operatorname{h}(R)$, to be the the number of hyperplanes in $\mathsf{Shi}^{(m)}(\Phi)$ that separate $R$ from $R^0$ and the \emph{coheight} by $\operatorname{coh}(R) := mN-\operatorname{h}(R)$. Furthermore, we (combinatorially) define $q$\emph{-Fu\ss -Catalan numbers} associated to $W$ by
		$$\operatorname{Cat}^{(m)}(W,q) := \sum_R{q^{\operatorname{coh}(R)}},$$
		where the sum ranges over all positive regions of $\mathsf{Shi}^{(m)}(\Phi)$.
	\end{definition}
	\begin{example}
		\begin{figure}
			\centering
			\includegraphics[width=170pt]{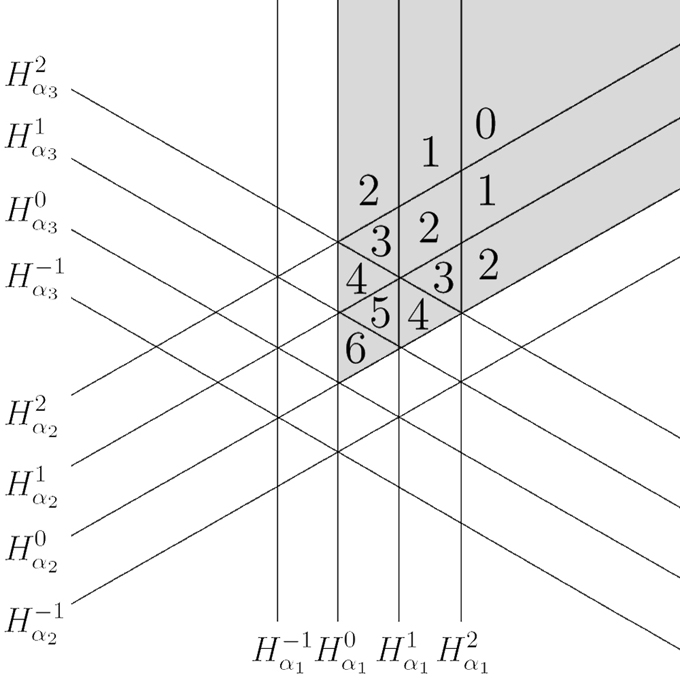}
			\caption{The extended Shi arrangement of type $A_2$ and $m=2$.}
			\label{arrangement}
		\end{figure}
		Let $W$ be the reflection group of type $A_2$ and $m=2$. In Fig.~\ref{arrangement}, the extended Shi arrangement of the given type is shown. The positive roots are denoted by $\alpha_1$, $\alpha_2$ and $\alpha_3=\alpha_1+\alpha_2$, the fundamental chamber is shaded and the positive regions are labelled by their coheights. This gives
		$$\operatorname{Cat}^{(2)}(W_{A_2},q) = 1+2q+3q^2+2q^3+2q^4+q^5+q^6.$$
	\end{example}
	\begin{remark} The poset of all regions, defined by the covering relation $R \prec R'$ if $\operatorname{h}(R)=\operatorname{h}(R')-1$ and $R, R'$ share a common face, is isomorphic to the poset $N\!N^{(k)}(W)$ described in \cite[Definition 5.1.20]{armstrong}.
	\end{remark}
	\begin{proposition}
	Let $I_1 \subseteq \ldots \subseteq I_m \subseteq \Phi^+$ be a filtered chain and let $R$ be the associated region. Then the bijection given above implies that
	$$\operatorname{coh}(R) = \sum{\# I_i}.$$
	\end{proposition}
	The next theorem shows that the coheight on regions in the fundamental chamber reduces for type $A$ to the area on $m$-Catalan paths of type $A$ and for type $B$ with $m=1$ to the area on Catalan paths of type $B$.
	\begin{theorem}\label{bij} For all integers $m \geq 1$, we have
		\begin{eqnarray*}
			\operatorname{Cat}^{(m)}(W_{A_{n-1}},q) = \operatorname{Cat}^{(m)}_n(q) &,&\operatorname{Cat}^{(1)}(W_{B_n},q) = \operatorname{Cat}_{B_n}(q).
		\end{eqnarray*}
	\end{theorem}
	\begin{remark}
		Counting lattice paths consisting of north and east steps having a boundary can always be seen as counting order ideals in very special kinds of posets. The posets occurring for $\operatorname{Cat}^{(m)}(W,q)$ with $W$ of type $D_n$ with $n \geq 4$ and of type $B_n$ with $n,m \geq 2$ fail to have this property.
	\end{remark}
	 The definition of \emph{coheight} is motivated by the following conjecture which would partially answer Open Problem \ref{problem1}:
	\begin{conjecture}\label{conjecture4}
		Let $W$ be a crystallographic reflection group. Then
		$$\operatorname{Cat}^{(m)}(W,q,1) = \operatorname{Cat}^{(m)}(W,q).$$
	\end{conjecture}
	\subsection{Non-crystallographic reflection groups}
		So far, no definition for root posets for non-crystallographic reflection groups is known. In \cite{armstrong}, D. Armstrong suggests, how these root posets \emph{should} look like in types $I_2(m)$ and $H_3$. Our computations confirm these ideas: Let $W$ be the a reflection group of one of the following types: $I_2(5),I_2(10),I_2(12),H_3$. Then
		$$\operatorname{Cat}^{(m)}(W,q,1) = \sum_{I \trianglelefteq \hspace{2pt} \Phi^+}{q^{\# I}},$$
		where $\Phi^+$ is Armstrong's suggested root poset of type $W$. The situation in the cyclic group of order $k$, $W=G(k,1,1)$ is the following: The $q,t$-Fu\ss -Catalan numbers associated to $G(k,1,1)$ are all equal. As we already know the classical case $k=2$, we get
		$$\operatorname{Cat}^{(m)}(W,q,t) = \operatorname{Cat}^{(m)}_2(q,t)=q+t.$$
		We also computed $\operatorname{Cat}^{(m)}(W,q,t)$ for $W$ the non-well-generated reflection group of type $G(4,2,2)$ and the result was, up to $m=3$,
		$$\operatorname{Cat}^{(m)}(W,q,t) = \operatorname{Cat}^{(m)}(W_{G_2},q,t).$$
	\section{Connections to rational Cherednik algebras}\label{sec4}
		Let $W$ be a real reflection group or, equivalently, let $W$ be a (finite) Coxeter group acting on the polynomial ring $\mathbb{C}[\mathbf{x},\mathbf{y}]$ and let $M^{(m)}$ be the $W$-module defined in Section \ref{defqtfusscat}. It is easy to see that $M^{(m)}$ is, except for the sign twist, equal to the alternating component of the module
		$$\mathbf{R}[\mathbf{x},\mathbf{y}] := (I^{m-1}/I^{m-1}J) \otimes \epsilon^{m-1},$$
		where $J$ is the ideal generated by all \emph{invariant polynomials without constant term},
		$$M^{(m)} \cong (\mathbf{R}[\mathbf{x},\mathbf{y}] \otimes \epsilon)^W = \mathbf{e}_\epsilon\mathbf{R}[\mathbf{x},\mathbf{y}].$$
		\vspace{-8pt}
		\begin{remark}
			$\mathbf{R}[\mathbf{x},\mathbf{y}]$ reduces for $m=1$ to the \emph{coinvariant ring} $\mathbb{C}[\mathbf{x},\mathbf{y}] / J$ and was, in type $A$, introduced in \cite{garsiahaiman}.
		\end{remark}
		The conjectured connection to \emph{rational Cherednik algebras} is the following:  in \cite{dunklopdam}, C.F. Dunkl and E. Opdam constructed a certain $W$-module $L$ depending on a non-negative integer $m$. This module carries a natural tensor product filtration and in \cite[Theorem 1.6]{BEG} Y. Berest, P. Etingof and V. Ginzburg showed that the Hilbert series of the trivial component of its associated graded module $\operatorname{gr}(L)$ is given by
		$$\mathcal{H}(\mathbf{e}(\operatorname{gr}(L));q) = q^{-mN} \prod_{i=1}^l{\frac{[d_i+mh]_q}{[d_i]_q}}.$$
		Using \cite[Lemma 6.7 (2)]{gordonstafford} together with equalities $(7.7)$ and $(7.8)$ in \cite{BEG}, we obtain the following result which partially generalizes \cite[Theorem 5]{gordon}.
		\begin{theorem}\label{theogordon}
			Let $W$ be a real reflection group and let $\mathbf{R}[\mathbf{x},\mathbf{y}]$ be graded by degree in $\mathbf{x}$ minus degree in $\mathbf{y}$. Then there exists a natural surjection of graded $W$-modules,
			$$\mathbf{R}[\mathbf{x},\mathbf{y}] \otimes \epsilon \twoheadrightarrow \operatorname{gr}(L).$$
		\end{theorem}
		\begin{conjecture}\label{conjecture3}
			The kernel of this surjection does not contain a copy of the trivial representation.
		\end{conjecture}
		This conjecture would imply the following corollary:
		\begin{corollary}\label{conjecture3corollary}
			Let $M^{(m)}$ be graded by degree in $\mathbf{x}$ minus degree in $\mathbf{y}$. Then
			$$M^{(m)} \cong \mathbf{e}(\operatorname{gr}(L))$$
			as graded $W$-modules.
		\end{corollary}
		\begin{remark}
			Together with the above discussion, Corollary \ref{conjecture3} would imply Conjecture \ref{conjecture2}.
		\end{remark}

\end{document}